\documentclass[10 pt, journal, final, twocolumn]{IEEEtran}
\usepackage{enumerate}
\usepackage{amsmath, amsthm, amssymb}
\usepackage{float}
\usepackage{subcaption}
\usepackage{flexisym}
\usepackage{mathtools}
\usepackage{breqn}
\usepackage{array}
\usepackage{tikz}
\usepackage{cite}
\usetikzlibrary{shapes.geometric, arrows}
\usepackage[font = small]{caption}
\usepackage{paralist}
\usepackage{mathrsfs}
\usepackage{algorithm}
\usepackage{algpseudocode}
\usepackage[inline]{enumitem}
\usepackage{bbm}
\allowdisplaybreaks[1]

\newtheorem{lemma}{Lemma}
\newtheorem{theorem}{Theorem}

\newtheorem{assumption}{Assumption}

\DeclareMathOperator*{\trace}{trace}
\DeclareMathOperator*{\diag}{diag}

\renewcommand{\IEEEQED}{\IEEEQEDclosed}

\begin{document}
\title{Resilient Distributed Estimation: Sensor Attacks}
\author{Yuan Chen, Soummya Kar, and Jos\'{e} M. F. Moura
\thanks{Yuan Chen {\{(412)-268-7103\}}, Soummya Kar {\{(412)-268-8962\}}, and Jos\'{e} M.F. Moura {\{(412)-268-6341, fax: (412)-268-3890\}} are with the Department of Electrical and Computer Engineering, Carnegie Mellon University, Pittsburgh, PA 15217 {\tt\small \{yuanche1, soummyak, moura\}@andrew.cmu.edu}}
\thanks{This material is based upon work supported by the Department of Energy under Award Number DE-OE0000779 and by DARPA under agreement numbers DARPA FA8750-12-2-0291 and DARPA HR00111320007. The U.S. Government is authorized to reproduce and distribute reprints for Governmental purposes notwithstanding any copyright notation thereon. The views and conclusions contained herein are those of the authors and should not be interpreted as necessarily representing the official policies or endorsements, either expressed or implied, of DARPA or the U.S. Government.}}
\maketitle

\begin{abstract}
	This paper studies multi-agent distributed estimation under sensor attacks. Individual agents make sensor measurements of an unknown parameter belonging to a compact set, and, at every time step, a fraction of the agents' sensor measurements may fall under attack and take arbitrary values. We present  the Saturated Innovation Update ($\mathcal{SIU}$) algorithm for distributed estimation resilient to sensor attacks. Under the iterative $\mathcal{SIU}$ algorithm, if less than one half of the agent sensors fall under attack, then, all of the agents' estimates converge at a polynomial rate (with respect to the number of iterations) to the true parameter. The resilience of $\mathcal{SIU}$ to sensor attacks does not depend on the topology of the inter-agent communication network, as long as it remains connected. We demonstrate the performance of $\mathcal{SIU}$ with numerical examples. 
\end{abstract}


\section{Introduction}\label{sect: intro}

This paper studies resilient distributed parameter estimation under sensor attacks. A network of agents makes measurements of an unknown parameter $\theta^*$, while an attacker corrupts a subset of the measurements. The agents' goal is to recover the parameter in a \textit{fully distributed} setting by using local sensor measurements and exchanging information with neighbors over a communication network. Our work addresses the following questions:
\begin{enumerate*}
	\item In a fully distributed setting, how can agents exchange and combine local information to resiliently estimate $\theta^*$ even when a fraction of the agents' sensors are under attack?
	\item What is the maximum tolerable fraction of attacked sensors under which a distributed algorithm resiliently estimates $\theta^*$?
\end{enumerate*}
To this end, we develop the Saturated Innovation Update ($\mathcal{SIU}$) algorithm, a \textit{consensus+innovations} type algorithm~\cite{Kar1, Kar2, Kar3, Kar4} for resilient distributed parameter estimation,  and we provide sufficient conditions, on the fraction of agents that may be attacked, under which the $\mathcal{SIU}$ algorithm ensures that \textit{every} agent correctly estimates $\theta^*$.


\subsection{Related Work}
The Byzantine Generals problem~\cite{ByzantineGenerals} demonstrates the effect of adversarial agents in multi-agent consensus, a special case of distributed inference, in an all-to-all communication setting. Existing work has studied Byzantine attacks (i.e., misbehaving agents) in the context of numerous decentralized (in which agents communicate measurements to a fusion center) problems, e.g., references~\cite{Varshney1, Varshney2, FusionCenterDetection2} study hypothesis testing, and reference~\cite{FusionCenterEstimation1} studies parameter estimation.  The resilience of the algorithms proposed in~\cite{Varshney1, Varshney2, FusionCenterEstimation1, FusionCenterDetection2} depends on the fraction of Byzantine agents. 


For fully distributed setups, in which there is no fusion center, the authors of~\cite{LeBlanc1, ResilientConsensus, LeBlanc2, Sundaram3} propose local filtering algorithms, in which each agent ignores, at every time step, a subset of its received messages, for consensus~\cite{LeBlanc1, ResilientConsensus}, scalar parameter estimation~\cite{LeBlanc2}, and optimization~\cite{Sundaram3} with Byzantine agents. The resilience of the algorithms proposed in~\cite{LeBlanc1, ResilientConsensus, LeBlanc2, Sundaram3} depends on the total number of adversarial agents and the topology of the inter-agent communication network. Currently, however, there is no computationally efficient way to evaluate and design resilient topologies~\cite{LeBlanc1}.  

The authors of~\cite{Pasqualetti1} design an algorithm for resilient consensus in which each agent needs to know the entire network topology. Similarly, reference~\cite{Sundaram1} proposes a method for resilient distributed  calculation of a function in which each agent uses its knowledge of the entire network topology to detect and identify adversaries. The resilience of the algorithms proposed in~\cite{Pasqualetti1} and~\cite{Sundaram1} depends on the total number of Byzantine agents as well as the connectivity of the inter-agent communication network. Reference~\cite{KalmanFailure} provides methods for sensor placement in response to sensor failures, i.e., when sensors do not produce any measurements, in distributed Kalman filtering.\footnote{In contrast, this paper addresses estimation with sensor attacks, i.e., when sensors produce arbitrary measurements as determined by an adversary.} Our previous work~\cite{ChenDistributed1} proposed a distributed algorithm for parameter estimation and detection of Byzantine agents that depends only on local knowledge and local communication. The resilience of the algorithm proposed in~\cite{ChenDistributed1} depends on the connectivity and global observability of the normally behaving agents.


\subsection{Summary of Contributions}
This paper presents the Saturated Innovation Update ($\mathcal{SIU}$) algorithm, an iterative, \textit{consensus+innovations} algorithm for resilient distributed parameter estimation. We consider a multi-agent setup, where each agent makes streams of measurements (over time) of a parameter $\theta^*$, and, at each time step an attacker manipulates a subset of the measurements. Each agent maintains a local estimate of $\theta^*$ and updates its estimate using a weighted combination of its neighbors' estimates (consensus) and its sensor measurement (innovation). 

During each time step (iteration), each agent applies a time-varying scalar gain to its own innovation term to ensure that the $\ell_2$ norm of the scaled innovation term is below a given threshold. The $\mathcal{SIU}$ algorithm only requires agents to have local knowledge of the network topology (i.e., agents only need to know their neighbors and not the entire network topology), and, achieves the same level of resilience as the most resilient centralized estimator. As long as less than half of the agents' sensors are under attack at any time step, then, all agents recover $\theta^*$ correctly.



The $\mathcal{SIU}$ algorithm offers two main advantages over existing techniques for resilient distributed computation under Byzantine attacks~\cite{Pasqualetti1, LeBlanc1, LeBlanc2, Sundaram1, Sundaram3, ResilientConsensus}. First, unlike the algorithms presented in~\cite{Pasqualetti1} and~\cite{Sundaram1}, the $\mathcal{SIU}$ algorithm only depends on agents having local knowledge and do not require agents to know the entire network. Second, under the $\mathcal{SIU}$ algorithms, the number of tolerable attacks at each time step scales linearly with the total number of agents regardless of the network topology, so long as the network is connected. To the best of our knowledge, there is no other algorithm, aside from $\mathcal{SIU}$, for distributed parameter estimation under sensor attack whose resilience does not depend on network topology. 


The rest of this paper is organized as follows. We review technical background and specify our sensing, communication, and attack models in Section~\ref{sect: background}. In Section~\ref{sect: algorithm}, we present the $\mathcal{SIU}$ algorithm. We state and prove our main results in Sections~\ref{sect: main} and~\ref{sect: MSIU}, respectively. Under $\mathcal{SIU}$, the agents' local estimates converge at a polynomial rate to $\theta^*$ when less than $\frac{1}{2}$ of the agents' sensors are under attack. We provide numerical examples of the $\mathcal{SIU}$ algorithm in Section~\ref{sect: example}, and we conclude in Section~\ref{sect: conclusion}.

\textit{Notation}: Let $\mathbb{R}^k$ be the $k$ dimensional Euclidean space, $I_k$ the $k$ by $k$ identity matrix, and $\mathbf{1}_k$ and $\mathbf{0}_k$ the column vectors of ones and zeros in $\mathbb{R}^k$, respectively. The operator $\left\lVert \cdot \right\rVert_2$ is the $\ell_2$ norm when applied to vectors and the induced $\ell_2$ norm when applied to matrices. For $v_1, v_2 \in \mathbb{R}^k$, let $\left\langle v_1, v_2 \right \rangle = v_1^T v_2$ be the inner product between $v_1$ and $v_2$. The Kronecker product of matrices $A$ and $B$ is $A \otimes B$ . For a symmetric matrix $\mathcal{M} = \mathcal{M}^T$, $\mathcal{M} \succeq 0$ ($M \succ 0$) means that $M$ is positive semi-definite (positive definite). 

Let $G = \left(V, E\right)$ be a simple, undirected graph, where $V = \left\{1, \dots, N\right\}$ is the set of vertices and $E$ is the set of edges. The neighborhood $\Omega_n$ of a vertex $n$ is the set of vertices that share an edge with $n$. Let $d_n = \left\vert \Omega_n \right\vert$ be the degree of a vertex, and the degree matrix of $G$ is $D = \diag\left( d_1, \dots, d_N \right)$. The adjacency matrix of $G$, $A = \left[ A_{nl} \right]$, where $A_{nl} = 1$ if $\left(n, l\right) \in E$ and $A_{nl} = 0$, otherwise, describes the structure of $G$. Let $L = D - A$ be the graph Laplacian of $G$. The eigenvalues of $L$ can be ordered as $0 = \lambda_1 (L) \leq \dots \leq \lambda_N(L)$, and $\mathbf{1}_N$ is the eigenvector associated with $\lambda_1(L)$. For a connected graph $G$, $\lambda_2 (L) > 0$. References~\cite{Spectral, ModernGraph} provide a detailed description of spectral graph theory.

\section{Background}\label{sect: background}

\subsection{Sensing and Communication Model}
Consider a network of $N$ agents $\left\{1, 2, \dots, N \right\}$ connected through a time-invariant inter-agent communication network $G = (V, E)$ (i.e., the topology of $G$ does not change over time). The vertex set $V = \left\{1, 2, \dots, N \right\}$ is the set of agents, and the edge set $E$ represents the inter-agent communication links. The agents' goal is to collectively estimate an unknown (non-random), static parameter $\theta^* \in \mathbb{R}^M$ (i.e., $\theta^*$ does not change over time). In the absence of an attacker, each agent $n$ makes a measurement $y_n(t)$ of the parameter $\theta^*$
\begin{equation}\label{eqn: sensingModel}
	y_n(t) = \theta^*,
\end{equation}
where $t$ is the (discrete) time index. Let $x_n(t)$ be the estimate of $\theta^*$ by agent $n$ at time $t$. The agents' goal is to ensure that 
\begin{equation}\label{eqn: agentGoal}
	x_n(t) \rightarrow \theta^*
\end{equation}
for all agents $n = 1, \dots, N$. Sensing model~\eqref{eqn: sensingModel} implies that, in the absence of an attacker, each agent is locally observable. That is, in the absence of an attacker, each agent, $n$, can exactly determine $\theta^*$ from $y_n(t)$. In the presence of an attacker, however, estimating $\theta^*$ becomes a nontrivial task.


We make the following assumptions regarding the sensing and network model.
\begin{assumption}\label{a: gConnected}
	The graph $G$ is connected.
\end{assumption}
\noindent Assumption~\ref{a: gConnected} can be made without loss of generality, since, if $G$ were not connected, we can separately consider each connected component of $G$. 
\noindent The focus of this paper is \textit{resilient} distributed estimation with respect to sensor data attacks, so, for simplicity and clarity, we assume that the inter-agent communication is noiseless. 
\begin{assumption}\label{a: thetaCompact}
	The $\ell_2$ norm of $\theta^*$ is bounded by a finite value $\eta$. That is, the parameter $\theta^*$ belongs to a set $\Theta$, defined as
	\begin{equation}\label{eqn: ThetaDef}
		\Theta = \left\{ \theta \in \mathbb{R}^M :  \left\lVert \theta \right \rVert_2 \leq \eta \right\}.
	\end{equation}
Each agent a priori knows the value of $\eta$. 
\end{assumption}
\noindent Boundedness is a natural assumption, since in many practical settings for distributed estimation, we estimate parameters from physical processes (e.g., power grid state estimation~\cite{Xie2} and wireless sensor networks for environmental monitoring~\cite{Kar2}\footnote{Existing work on distributed estimation (e.g.~\cite{Kar2, Xie2}) does not account for the effect of sensor attacks. In contrast, this paper provides a distributed estimation algorithm that is resilient to sensor attacks.}). The values of such parameters are bounded by laws of physics.

\subsection{Attacker Model}\label{sect: attackerModel}
An attacker aims to disrupt the estimation procedure and prevent the agents from achieving their goal~\eqref{eqn: agentGoal}. The attacker may arbitrarily manipulate the sensor measurements of some of the agents. We model the measurement of a sensor under attack as
\begin{equation}\label{eqn: attackedSensor}
	y_n(t) = \theta^* + a_n(t),
\end{equation}
where $a_n(t)$ is the disturbance induced by the attacker. In practice, an attacker does not directly design the additive disturbance $a_n(t)$. Instead, the attacker replaces the agents' measurement with any arbitrary $y_n(t)$, which can be modeled, following equation~\eqref{eqn: attackedSensor}, by a corresponding $a_n(t)$. The attacker may know the true value of the parameter $\theta^*$ and uses this information to determine how to manipulate sensor measurements and inflict maximum damage.\footnote{In general, an attacker does not \textit{need} to know the value of $\theta^*$ to manipulate measurements. An attacker who does not know the value of $\theta^*$ may not be able to inflict as much damage as one who does but, without proper countermeasures, can still prevent the agents from recovering $\theta^*$.} 
{\color{black} In our attacker model~\eqref{eqn: attackedSensor}, the adversary directly manipulates the sensor measurements of a subset of the agents. Since the adversary directly manipulates the measurements, model~\eqref{eqn: attackedSensor} is a type of \textit{spoofing} attack~\cite{FusionCenterEstimation1, Zhang2, Zhang3}. This model is related to man-in-the-middle attacks, where the adversary manipulates communications between the agents and sensors or manipulates data after processing (e.g., quanitization)~\cite{Zhang2, Zhang3}.} Let the set $\mathcal{A}_t$ denote the set of agents whose sensors are under attack at time $t$, i.e.,
\begin{equation}\label{eqn: calADef}
	\mathcal{A}_t = \left\{n \in V \vert a_n(t) \neq 0 \right\}.
\end{equation}
Let $\mathcal{N}_t = V\setminus \mathcal{A}_t$ denote the agents whose sensors are \textit{not} under attack at time $t$. In the presence of an attacker, agents can no longer achieve~\eqref{eqn: agentGoal} by setting $x_n(t) = y_n(t)$. 

We make the following assumptions on {\color{black}the attacker}:
\begin{assumption}\label{a: attackerAssumptions}
	For some $0 \leq S < N$, $\left \lvert \mathcal{A}_t \right \vert \leq S$, for all $t$.
\end{assumption}
{\color{black}\begin{assumption}\label{a: staticParam}
	The attack does not change the value of $\theta^*$. 
\end{assumption}}
\noindent The agents do not know the set $\mathcal{A}_t$. Assumption~\ref{a: attackerAssumptions} is similar to the sparse sensor attack assumption found in the cyber-physical security literature~\cite{ChenICASSP, Fawzi, Shoukry}. We will compare the resilience of the $\mathcal{SIU}$ algorithm to the resilience of any centralized estimation algorithm.  Applying {\color{black}Theorem 3.2 from~\cite{Shoukry}}, a necessary and sufficient condition for any centralized estimator to resiliently recover $\theta^*$ is that the number of attacked sensors is less than half of all sensors, i.e., $\frac{\left\lvert \mathcal{A}_t \right \rvert}{N} < \frac{1}{2}$. {\color{black} Assumption~\ref{a: staticParam} states that the attacker does not change the value of the parameter of interest. The attacker may only manipulate the measurements $y_n(t)$ of a subset of the agents.}


The attacker model~\eqref{eqn: attackedSensor} differs from the Byzantine Attacker Model~\cite{ByzantineGenerals}. A Byzantine attacker is able to hijack a subset of the agents and control \textit{all} aspects of their behavior (i.e., a Byzantine attacker can arbitrarily manipulate the message generation and estimate generation processes of hijacked agents). Equation~\eqref{eqn: attackedSensor} models a simpler attack than the Byzantine model, but it is still a realistic model. Due to resource limitations (e.g., time, computation power, etc., we refer the reader to~\cite{TeixeiraModels} for a detailed description relating attacker resources to attacker capabilities), the attacker directly manipulates sensor measurements but does not completely hijack individual agents. 

\section{Saturated Innovation Update ($\mathcal{SIU}$) Algorithm}\label{sect: algorithm}
The Saturated Innovation Update ($\mathcal{SIU}$) algorithm is a \textit{consensus+innovations}~\cite{Kar1, Kar2, Kar3, Kar4} algorithm for resilient distributed estimation. In $\mathcal{SIU}$, every agent maintains and updates a local estimate $x_n(t)$. In contrast with~\cite{Kar1, Kar2, Kar3, Kar4}, which assumes that there is no attacker, the $\mathcal{SIU}$ algorithm addresses distributed estimation when a subset of the agents' sensors are under attack.  Each iteration of $\mathcal{SIU}$ consists of two steps: message passing and estimate update. To initialize, each agent $n$ sets its local estimate $x_n(0) = 0$.

\underline{Message Passing}: In each iteration, every agent $n$ transmits its current estimate $x_n(t)$ to each of its neighbors. Each agent $n$ transmits $d_n$ messages in each iteration.

\underline{Estimate Update}: Each agent $n \in V$ updates its estimate as
\begin{equation}\label{eqn: agentUpdate}
\begin{split}
	x_n(t+1) =& x_n(t) - \beta_t \sum_{l \in \Omega_n} \left( x_n(t) - x_l(t)\right) \\ &+ \alpha_t K_n(t)\left(y_n - x_n(t) \right),
\end{split}
\end{equation} 
where $\alpha_t >0, \beta_t>0$ are sequences of parameters to be specified in the sequel. The term $K_n(t)$ is a time-varying gain on the local innovation and is defined as
\begin{equation}\label{eqn: saturationGain}
	K_n(t) = \left\{ \begin{array}{cl}1, & \left\lVert y_n(t) - x_n(t) \right \rVert_2 \leq \gamma_t \\
 \frac{\gamma_t}{\left\lVert y_n(t) - x_n(t) \right \rVert_2}, & \text{otherwise}, \end{array} \right.
\end{equation}
where $\gamma_t$ is a parameter to be specified in the sequel. 

The gain $K_n(t)$ ensures that the $\ell_2$ norm of the term $K_n(t) \left( y_n(t) - x_n(t) \right)$ is upper bounded by $\gamma_t$ for all $n \in V$ and for all $t$. The challenge in designing the $\mathcal{SIU}$ algorithm is selecting the adaptive threshold, $\gamma_t$, which represents how much sensor measurements are allowed to deviate from the local estimates. If $\gamma_t$ is too small, then, the gain $K_n(t)$ will limit the impact of uncompromised sensors ($n \in \mathcal{N}_t$), and the agents will not recover $\theta^*$. On the other hand, if $\gamma_t$ is too large, then, the attacker will be able to mislead the agents with measurements the deviate more from the local estimates. The key is to choose $\gamma_t$ to balance these two effects.

We adopt the following parameter selection procedure:
\begin{enumerate}
	\item Select resilience index $0<s<\frac{1}{2}$.

	\item Select the sequence $\alpha_t$ and $\beta_t$ to be of the form
	\begin{equation}\label{eqn: mtsAlphaBeta}
	\alpha_t = \frac{a}{(t+1)^{\tau_1}}, \: \beta_t =\frac{b}{(t+1)^{\tau_2}},
\end{equation}
where $ 0 < a \leq \frac{1}{1-2s}$, $0 < b \leq \frac{1}{\lambda_N\left(L\right)}$, $0 < \tau_2 < \tau_1 <1$.

\item \noindent The corresponding $\gamma_t$ sequence is
	\begin{equation}\label{eqn: mtsGamma}
		\gamma_t = \gamma_{1, t} + \gamma_{2, t},
	\end{equation}
	where $\gamma_{1, t}$ and $\gamma_{2, t}$ follow
	\begin{equation}\label{eqn: mtsGammaDynamics}
		\begin{split}
			\gamma_{1, t+1} &= \left(1 - \beta_t \lambda_2 \left(L\right) + \kappa_1 \alpha_t \right) \gamma_{1, t} + \kappa_2\alpha_t\gamma_{2, t},\\
			\gamma_{2, t+1} &=  \alpha_t \gamma_{1,t}+ \left(1 - \alpha_t \left(1-2s\right)\right) \gamma_{2, t},
		\end{split}
	\end{equation}

for $\kappa_1 = 1+\sqrt{N}$ and $\kappa_2 = {\color{black}2}\sqrt{N}$ with initial conditions $\gamma_{1, 0} = 0$ and $\gamma_{2, 0} = \eta$. Recall, from~\eqref{eqn: ThetaDef}, that $\eta$ is the upper bound on $\left\lVert \theta^*\right\rVert_2$. 
\end{enumerate}

\noindent The resilience index $s$ determines the maximum tolerable fraction of attacked sensors, and, for the $\mathcal{SIU}$ algorithm, we require $s < \frac{1}{2}$. The resilience index is a design parameter in $\mathcal{SIU}$; choosing a larger $s$ improves the resilience of $\mathcal{SIU}$ but slows the convergence of estimates to the true parameter.

The gains $K_{n}(t)$ are state-dependent (i.e., they depend on the estimates $x_n(t)$) and non-smooth (as functions of the state), which makes the analysis of $\mathcal{SIU}$ nontrivial and quite different from traditional consensus+innovations procedures~\cite{Kar1, Kar2, Kar3, Kar4} that rely on the smoothness of the gains to obtain appropriate Lyapunov conditions for convergence analysis. The analysis with non-smooth gains as in~\eqref{eqn: saturationGain} require new technical machinery that we develop in this paper.

\section{Main Result}\label{sect: main}
We now present our main result, which addresses the resilience and performance of the $\mathcal{SIU}$ algorithm.

\begin{theorem}[Resilience of Algorithm $\mathcal{SIU}$]\label{thm: mtsMain}
	Let $\alpha_t$, $\beta_t$, and $\gamma_t$ be given by~\eqref{eqn: mtsAlphaBeta},~\eqref{eqn: mtsGamma}, and~\eqref{eqn: mtsGammaDynamics}, and let $s \in \left(0, \frac{1}{2} \right)$ be the resilience index. Under the $\mathcal{SIU}$ algorithm, if $\frac{\left\lvert \mathcal{A}_t \right\rvert}{N} < s$ for all times $t$, then, we have
\begin{equation}\label{eqn: mtsMain}
	\lim_{t\rightarrow\infty} (t+1)^{\tau_0} \left\lVert x_n(t) - \theta^*\right\rVert_2 = 0,
\end{equation}
for every agent $n$ and for all $0 \leq \tau_0 < \tau_1 - \tau_2$.
\end{theorem}

\noindent Theorem~\ref{thm: mtsMain} states that, under $\mathcal{SIU}$, all of the agents' estimates converge to $\theta^*$ at a rate of ${1 \over t^{\tau_0}}$, for any $0 \leq \tau_0 < \tau_1 - \tau_2$, so long as less than half of the agents' sensors fall under attack, irrespective of \textit{how} the attacker manipulates the sensor measurements. Recall that any centralized estimator, which collects sensor measurements from all of the nodes at once, is only resilient to attacks on less than half of the sensors. {\color{black} The $\mathcal{SIU}$ algorithm ensures that all of the agents' recover the value $\theta^*$ even though it does not explicity identify those agents who are under attack. }The $\mathcal{SIU}$ algorithm achieves, in a fully distributed setting, the same resilience as the most resilient centralized estimator, regardless of the topology of the inter-agent communication network (so long as it is connected). 
%


Under the $\mathcal{SIU}$ algorithm, the message generation and estimate update for each agent only depends on the agent's local knowledge. Existing algorithms for resilient distributed consensus~\cite{Pasqualetti1} and resilient distributed function calculation~\cite{Sundaram1} require each agent to have knowledge of the entire structure of the communication network. When there are many agents in the network, algorithms requiring global network knowledge become expensive from both {\color{black} memory and computation perspectives}. 

In addition, existing work on distributed algorithms resilient to Byzantine attackers shows that, in general, the resilience of an algorithm (i.e., the number of tolerable Byzantine agents) depends on the structure of the network~\cite{Sundaram1, LeBlanc1, ResilientConsensus, LeBlanc2, Pasqualetti1, Sundaram3}. Under our attack model~\eqref{eqn: attackedSensor}, the number of tolerable attacks for $\mathcal{SIU}$ algorithms scales linearly with the number of agents, regardless of the (connected) network topology. When the resilience of algorithms depends on network structure (e.g.~\cite{LeBlanc1, LeBlanc2, Sundaram1, Sundaram3, ResilientConsensus, Pasqualetti1}), the number of tolerable attacks does not necessarily increase as the number of agents grows. In comparison, for \textit{any} connected structure, $\mathcal{SIU}$ ensures resilient estimation if less half of the agents' sensors are under attack.

\section{Analysis of $\mathcal{SIU}$: Proof of Theorem~\ref{thm: mtsMain}}\label{sect: MSIU}


In order to prove Theorem~\ref{thm: mtsMain}, we need to analyze convergence properties of time-varying linear systems of the form:
\begin{equation}\label{eqn: miStackedSystem}
	\begin{split}
		v_{t+1} &= \left(1 - c_3 r_1(t) \right) v_t + c_4 r_1(t) w_t, \\
		w_{t+1} &= \left(1 - c_5 r_2(t) + c_6 r_1(t) \right) w_t + c_7 r_1(t) v_t,
	\end{split}
\end{equation}
with nonzero initial conditions $v_0, w_0 \neq 0$ where $r_1(t)$ and $r_2(t)$ follow\begin{equation}\label{eqn: mtsSystem1Rates}
	r_1(t) = \frac{c_1}{(t+1)^{\delta_1}},\: r_2(t) = \frac{c_2}{(t+1)^{\delta_2}},
\end{equation} with $\delta_1 > \delta_2$ and $c_3, \dots, c_7 > 0$. Specifically, we require the following lemma, the proof of which may be found in the appendix.


\begin{lemma}\label{lem: mi4}
	The system in~\eqref{eqn: miStackedSystem} satisfies
	\begin{align}
		\lim_{t \rightarrow \infty} (t+1)^{\delta_0} v_t &= 0, \label{eqn: mi4}\\
		\lim_{t \rightarrow \infty} (t+1)^{\delta_0} w_t &= 0, \label{eqn: mi4a}
	\end{align}
	for all $0 \leq \delta_0 < \delta_1 - \delta_2$.
\end{lemma}

\noindent We now prove Theorem~\ref{thm: mtsMain}. 
\begin{IEEEproof}[Proof (Theorem~\ref{thm: mtsMain}):] We analyze the behavior of the agents' average estimate. Define the stacked estimate $\mathbf{x}_t = \left[\begin{array}{ccc} x_1(t)^T & \cdots & x_N(t)^T \end{array}\right]^T,$ the stacked measurements $\mathbf{y}_t = \left[\begin{array}{ccc} y_1(t)^T & \cdots & y_N(t)^T \end{array}\right]^T,$ and the network average estimate $\overline{\mathbf{x}}_t = \frac{1}{N} \left(\mathbf{1}_N^T \otimes I_M\right) \mathbf{x}_t.$ 	Further define 
	\begin{align}\label{eqn: ktilde}
		\widetilde{K}_n (t) &= K_n(t) \mathbbm{1}_{\left\{ n \in \mathcal{N}_t \right\}}, \\
		\mathbf{K}_t &= \diag\left(K_1(t), \dots, K_N(t) \right), \\
		\mathbf{K}_{\mathcal{N}_t} &=  \left(\widetilde{K}_1(t), \dots, \widetilde{K}_N (t) \right), \\
		\mathbf{K}_{\mathcal{A}_t} &= \mathbf{K}_t - \mathbf{K}_{\mathcal{N}_t}.
	\end{align}


\underline{Step 1}: We determine the dynamics of $\mathbf{x}_t$ and $\overline{\mathbf{x}}_t$. From~\eqref{eqn: agentUpdate}, we have that $\mathbf{x}_t$ follows the dynamics
\begin{equation}\label{eqn: mtsProof1}
\begin{split}
	\mathbf{x}_{t+1} &= \left( I_{NM} - \beta_t L \otimes I_M  - \alpha_t \mathbf{K}_{\mathcal{N}_t} \otimes I_M \right) \mathbf{x}_t + \\
	&\quad \alpha_t \left(\mathbf{K}_{\mathcal{A}_t} \otimes I_M \right) \mathbf{y}_t.
\end{split}
\end{equation}
Define $\widehat{\mathbf{y}}_t$, the difference between the local estimates and the average estimate, as
	$\widehat{\mathbf{y}}_t = \mathbf{x}_t - \mathbf{1}_N \otimes \overline{\mathbf{x}}_t,$
and $\overline{\mathbf{e}}_t$, the average error, as $\overline{\mathbf{e}}_t = \overline{\mathbf{x}}_t - \theta^*$.
Further define the matrix $P_{NM}$ as
\begin{equation}\label{eqn: pnmDef}
	P_{NM} = \frac{1}{N} \left( \mathbf{1}_N \otimes I_M \right) \left( \mathbf{1}_N \otimes I_M \right)^T.
\end{equation}
Note that
	$P_{NM} \left(\mathbf{1}_N \otimes \overline{\mathbf{x}}_t \right) = \mathbf{1}_N \otimes \overline{\mathbf{x}}_t.$
From~\eqref{eqn: mtsProof1}, we have that $\widehat{\mathbf{y}}_t$ follows
\begin{equation}\label{eqn: mtsProof6}
\begin{split}
	&\widehat{\mathbf{y}}_{t+1} = \left(I_{NM} - \beta_t L \otimes I_M - P_{NM} \right)\widehat{\mathbf{y}}_t + \\
	&\alpha_t \left(I_{NM} - P_{NM} \right) \left(\left( \mathbf{K}_{\mathcal{N}_t} + \mathbf{K}_{\mathcal{A}_t} \right) \otimes I_M\right) \left( \mathbf{y}_t - \mathbf{x}_t \right).
\end{split}
\end{equation}
Following~\cite{Kar3}, for $0 < \beta_t \leq \frac{1}{\lambda_N (L)}$, the eigenvalues of the matrix $I_{NM} - \beta_t L \otimes I_M - P_{NM}$ are $0$ and $1 - \beta_t \lambda_n \left(L\right)$ for $n = 2, \dots, N$, each repeated $M$ times.

For $n \in \mathcal{N}_t$, we have $\mathbf{y}_t - \mathbf{x}_t = \widehat{\mathbf{y}}_t + \mathbf{1}_N \otimes \overline{\mathbf{e}}_t$. Then, equation~\eqref{eqn: mtsProof6} becomes
\begin{equation}\label{eqn: mtsProof7}
\begin{split}
	&\widehat{\mathbf{y}}_{t+1} = \left(I_{NM} - \beta_t L \otimes I_M - P_{NM} \right)\widehat{\mathbf{y}}_t - \\
	&\quad \alpha_t \left(I_{NM} - P_{NM} \right)\left( \mathbf{K}_{\mathcal{N}_t} \otimes I_M \right)\left( \widehat{\mathbf{y}}_t + \mathbf{1}_N \otimes \overline{\mathbf{e}}_t \right) + \\
&\quad \alpha_t \left(I_{NM} - P_{NM} \right) \left(\mathbf{K}_{\mathcal{A}_t} \otimes I_M \right) \left( \mathbf{y}_t - \mathbf{x}_t \right).
\end{split}
\end{equation}
The process $\overline{\mathbf{e}}_t$ evolves as
\begin{equation}\label{eqn: mtsProof8}
\begin{split}
	\overline{\mathbf{e}}_{t+1} &= \left(1 - \frac{\alpha_t}{N}\trace\left(\mathbf{K}_{\mathcal{N}_t} \right) \right) \overline{\mathbf{e}}_t  - \\
&\quad \frac{\alpha_t}{N}\left( \mathbf{1}_N^T \otimes I_M \right) \left(\mathbf{K}_{\mathcal{N}_t} \otimes I_M \right) \widehat{\mathbf{y}}_t +\\
&\quad \frac{\alpha_t}{N}\left( \mathbf{1}_N^T \otimes I_M \right) \left(\mathbf{K}_{\mathcal{A}_t} \otimes I_M \right) \left( \mathbf{y}_t - \mathbf{x}_t \right).
\end{split}
\end{equation}

\underline{Step 2}: From~\eqref{eqn: mtsProof7} and~\eqref{eqn: mtsProof8}, we find upper bounds on $\left\lVert \widehat{\mathbf{y}}_t \right\rVert_2$ and $\left \lVert \overline{\mathbf{e}}_t \right \rVert_2$. Let $V_t = \left\lVert \widehat{\mathbf{y}}_t \right\rVert_2$, and let $W_t = \left\lVert \overline{\mathbf{e}}_t \right\rVert_2$. In order to compute $V_t$ and $W_t$, we require the following facts:
\begin{align}
	\left\lVert I_{NM} - P_{NM} \right \rVert_2 &= 1,\label{eqn: mpi1}\\
	\left\lVert \mathbf{K}_{\mathcal{N}_t} \otimes I_M \right \rVert_2 &\leq 1, \label{eqn: mpi2}\\
	\left\lVert \left(\mathbf{K}_{\mathcal{N}_t} \otimes I_M \right) \left( \mathbf{1}_N \otimes \overline{\mathbf{e}}_t \right) \right\rVert_2 & \leq \sqrt{\left\lvert \mathcal{N}_t \right \rvert} \left\lVert \overline{\mathbf{e}}_t \right \rVert_2. \label{eqn: mpi3}
\end{align}
Equation~\eqref{eqn: mpi1} follows directly from the definition of $P_{NM}$. Inequality~\eqref{eqn: mpi2} follows from the fact that $0 < K_n(t) \leq 1$. Inequality~\eqref{eqn: mpi3} follows from
\begin{align}
	\left\lVert \left(\mathbf{K}_{\mathcal{N}_t} \otimes I_M \right) \left( \mathbf{1}_N \otimes \overline{\mathbf{e}}_t \right) \right\rVert_2 &= \sqrt{\sum_{n \in \mathcal{N}_t} \left(K_n(t)\right)^2 \left\lVert \overline{\mathbf{e}}_t \right \rVert_2}, \nonumber \\
	& \leq \sqrt{\left\lvert \mathcal{N}_t \right \rvert} \left\lVert \overline{\mathbf{e}}_t \right \rVert_2.
\end{align}
Computing the $\ell_2$ norm of both sides of~\eqref{eqn: mtsProof7} and~\eqref{eqn: mtsProof8}, and using the fact that, $\left\lVert K_n(t) \left( y_n(t) - x_n(t)\right) \right \rVert_2 \leq \gamma_t = \gamma_{1, t} + \gamma_{2, t}$ for any $n \in \mathcal{A}_t$, we have
\begin{align}
\begin{split}\label{eqn: mtsProof12}
	V_{t+1} &\leq \left(1 - \beta_t \lambda_2 \left(L \right) + \alpha_t \right) V_t + \alpha_t W_t \sqrt{\left\lvert \mathcal{N}_t \right \rvert} + \\
&\quad \alpha_t \gamma_t\sqrt{\left\lvert \mathcal{A}_t \right \rvert},
\end{split}\\
\begin{split}\label{eqn: mtsProof13}
	W_{t+1} & \leq \left(1 - \frac{\alpha_t}{N}\trace\left(\mathbf{K}_{\mathcal{N}_t} \right) \right) W_t + \frac{\alpha_t \gamma_t \left\lvert \mathcal{A}_t \right \rvert}{N} + \\
&\quad  \frac{\alpha_t}{N}\left\lVert\left( \mathbf{1}_N^T \otimes I_M \right) \left(\mathbf{K}_{\mathcal{N}_t} \otimes I_M \right) \widehat{\mathbf{y}}_t \right \rVert_2.
\end{split}
\end{align}

\underline{Step 3}: We now use induction to show that, if $\frac{\left\lvert \mathcal{A}_t \right \rvert}{N} < s$ for all $t \geq 0$,  then, $V_t \leq {\gamma_{1, t}}$ and $W_t \leq {\gamma_{2, t}}$ for all $t \geq 0$. In the base case, we consider $t = 0$, and we have $x_n(0) = 0$ for all $n$. Thus, we have $\overline{\mathbf{x}}_0 = 0$. For all $n \in V$, $\left\lVert x_n(0) - \overline{\mathbf{x}}_0 \right\Vert_2 = 0$, which means that $\left\lVert \widehat{\mathbf{y}}_t \right\rVert_2 = 0 = \gamma_{1, 0}$. Moreover, since $\overline{\mathbf{x}}_0 = 0$, we have $\left\lVert \overline{\mathbf{e}}_t \right\rVert_2 = \left\lVert \theta^* \right\rVert_2 \leq \eta = \gamma_{2, 0}$.

\underline{Step 4}: In the induction step, we assume that $V_{t} \leq \gamma_{1, t}$ and $W_t \leq \gamma_{2, t}$, and we show that $V_{t+1} \leq \gamma_{1, t+1}$ and $W_{t+1} \leq \gamma_{2, {t+1}}$. Substituting the induction hypotheses into~\eqref{eqn: mtsProof12} and~\eqref{eqn: mtsProof13}, we have
\begin{align}
\begin{split}\label{eqn: mtsProof16}
	V_{t+1} &\leq \left(1 - \beta_t \lambda_2 \left(L \right) + \kappa_1\alpha_t \right)\gamma_{1, t} + 2\alpha_t \gamma_{2, t} \sqrt{N}, 
\end{split}\\
\begin{split}\label{eqn: mtsProof17}
	W_{t+1} &\leq \left(1 - \alpha_t \left(\trace\left(\frac{\mathbf{K}_{\mathcal{N}_t}}{N}\right)  - \frac{\left\lvert\mathcal{A}_t\right\rvert}{N}\right)\right)\gamma_{2, t} + \\
	& \quad \alpha_t \left(\frac{\left\lvert \mathcal{N}_t \right \rvert + \left\lvert \mathcal{A}_t\right\rvert}{N} \right) \gamma_{1, t},
\end{split}
\end{align}
where, recall, $\kappa_1 = 1 + \sqrt{N}$. 
To derive~\eqref{eqn: mtsProof16}, we have used the fact that ${\color{black}2}\sqrt{N} \geq \sqrt{\left\lvert \mathcal{A}_t \right \rvert} + \sqrt{\left\lvert \mathcal{N}_t \right \rvert}$. To derive~\eqref{eqn: mtsProof17}, we have used the fact that $\left\lVert\left( \mathbf{1}_N^T \otimes I_M \right) \left(\mathbf{K}_{\mathcal{N}_t} \otimes I_M \right) \widehat{\mathbf{y}}_t \right \rVert_2 \leq \left\lvert \mathcal{N}_t \right \rvert \left\lVert \widehat{\mathbf{y}}_t \right \rVert_2.$ Inequality~\eqref{eqn: mtsProof17} states that $W_{t+1}$ depends on the term \[\trace\left(\frac{\mathbf{K}_{\mathcal{N}_t}}{N} \right) = \frac{1}{N} \sum_{n \in \mathcal{N}_t} K_n(t),\]
which, by definition of $K_n(t)$, depends on $\left\lVert y_n(t) - x_n(t) \right \rVert_2$. Note that, by definition of $\mathcal{N}_t$, we have $y_n(t) = \theta^*$ for all $n \in \mathcal{N}_t$. Applying the triangle inequality and the induction hypotheses, we have
\begin{align}
	\left\lVert x_n(t) - \theta^* \right \rVert_2 & \leq \left\lVert x_n(t) - \overline{\mathbf{x}}_t \right \rVert_2 + \left \lVert \overline{\mathbf{x}}_t - \theta^* \right \rVert_2, \label{eqn: mtsProof18}\\
	& \leq \gamma_{1, t} + \gamma_{2, t} = \gamma_t. \label{eqn: mtsProof19}
\end{align}
Inequality~\eqref{eqn: mtsProof19} implies that $K_n(t) = 1$ for all $n \in \mathcal{N}_t$, which means that
	$\trace\left(\frac{\mathbf{K}_{\mathcal{N}_t}}{N} \right) = 1 - \frac{\left\lvert \mathcal{A}_t \right \rvert}{N}.$
Substituting $\trace\left(\frac{\mathbf{K}_{\mathcal{N}_t}}{N} \right) = 1 - \frac{\left\lvert \mathcal{A}_t \right \rvert}{N}$ into~\eqref{eqn: mtsProof17} and using the fact that $\frac{\left\lvert\mathcal{A}_t \right\rvert}{N} < s$, we have
\begin{equation}\label{eqn: mtsProof21}
	W_{t+1} \leq \left(1 - \alpha_t \left(1-2s\right)\right) \gamma_{2, t} + \alpha_t \gamma_{1,t},
\end{equation}
which yields the relation $W_{t+1} \leq \gamma_{2, t}$. The relation $V_{t+1} \leq \gamma_{1, t}$ follows directly from~\eqref{eqn: mtsProof16}.

\underline{Step 5}: We now study the behavior of $\gamma_{1, t}$ and $\gamma_{2, t}$. So long as $0 < s < \frac{1}{2}$, system~\eqref{eqn: mtsGammaDynamics}, which describes the dynamics of $\gamma_{1, t}$ and $\gamma_{2, t}$, falls under the purview of Lemma~\ref{lem: mi4}.\footnote{If $s \geq {1 \over 2}$, then the right hand side of $\eqref{eqn: mtsProof21}$ does not converge to $0$. Thus,~\eqref{eqn: mtsProof21} establishes the maximum fraction of tolerable attacks for $\mathcal{SIU}$.} Thus we have
\begin{align}
	\lim_{t\rightarrow \infty} (t+1)^{\tau_0} \gamma_{i, t} &= 0, \label{eqn: mtsProof22}
\end{align}
$i = 1, 2,$ for all $0 \leq \tau_0 < \tau_1 - \tau_2$. Combining~\eqref{eqn: mtsProof22} with~\eqref{eqn: mtsProof19} yields the desired result: for every agent $n$ and for all $0 \leq \tau_0 < \tau_1 - \tau_2$,
\begin{equation}\label{eqn: mtsProof23}
	\lim_{t\rightarrow\infty} (t+1)^{\tau_0} \left\lVert x_n(t) - \theta^*\right\rVert_2 = 0.
\end{equation}
\end{IEEEproof}

As a consequence of Theorem~\ref{thm: mtsMain}, if $\frac{\left\lVert \mathcal{A}_t \right \rVert}{N} < s < \frac{1}{2}$ for all $t$, then, for all $n \in V$ and for all $0 \leq \tau_0 < \tau_1 - \tau_2$, there exists finite $R_n > 0$ such that 
	$\left\lVert x_n(t) - \theta^* \right \rVert_2 < \frac{R_n}{(t+1)^{\tau_0}}.$
The rate $\tau_0$ depends only on the choice of $\tau_1$ and $\tau_2$, but the constant $R_n$ depends on the behavior of $\gamma_{1, t}$ and $\gamma_{2, t}$. The behavior of $\gamma_{2, t}$ depends on the resilience index $s$ \textendash~the term $\left(1 - \alpha_t(1-2s)\right)$ increases as $s$ increases, which means that the constant $R_n$ increases as well. Thus, in the $\mathcal{SIU}$ algorithm, there is a trade off between resilience and convergence of estimates to the true parameter.

\section{Numerical Simulations}\label{sect: example}
\begin{figure}[h!]
\centering
\includegraphics[width = 0.5\columnwidth]{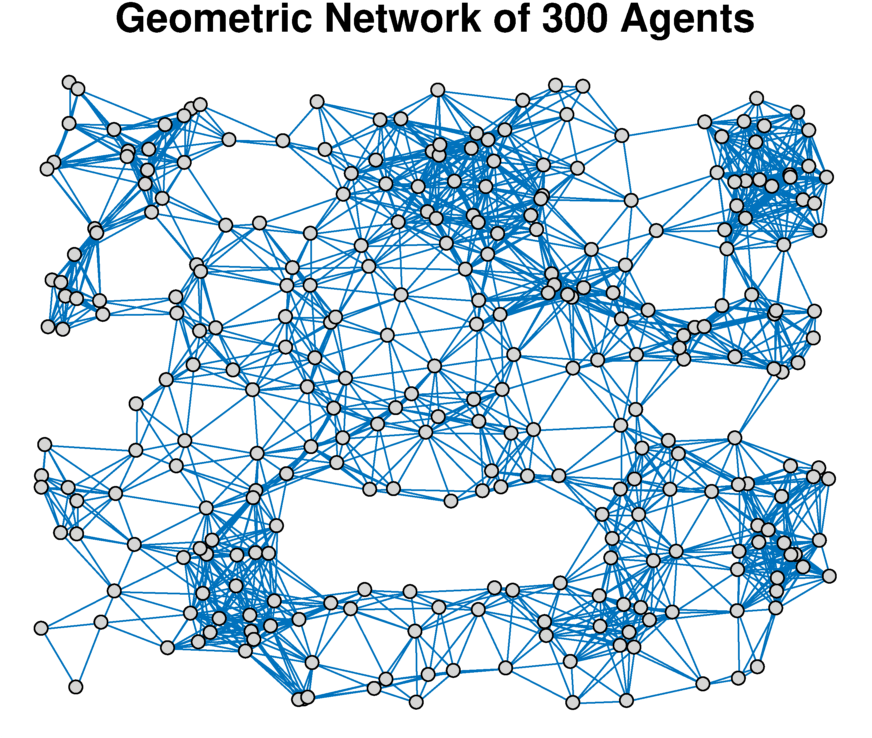}
\caption{Random geometric network of $N=300$ agents.}\label{fig: network}
\end{figure}

In our numerical simulations, we consider the random geometric network of $N = 300$ agents given by Figure~\ref{fig: network}. Each agent observes a parameter $\theta^* \in \mathbb{R}^{3}$ with bounded energy $\left\lVert \theta^* \right \rVert_2 \leq \eta = 100$. The agents represent, for example, a team of robots tracking a target. The parameter $\theta^*$ is the target's location, expressed in its $x$, $y$, and $z$ coordinates. We consider attacks on $S_1 = 60$ and $S_2 = 120$ agents, respectively, and we choose corresponding resilience indices of $s_1 = 0.201$ and $s_2 = 0.401$. In each simulation, we run $\mathcal{SIU}$ for $t = 500,000$ iterations. We choose the following parameters: $a = 1.54 \times 10^{-4}, b = 3.78 \times 10^{-2}, \tau_1 = 0.15, \tau_2 = 0.001$.

\begin{figure}[h!]
	\centering
	\includegraphics[width = \columnwidth]{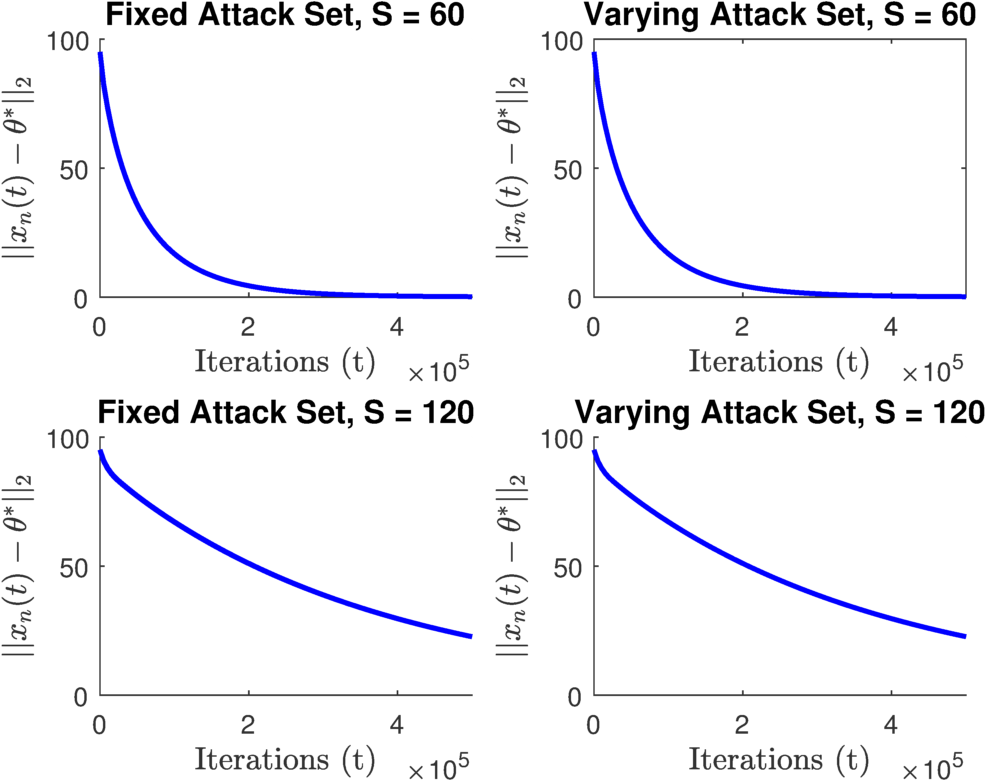}
	\caption{Performance of $\mathcal{MSIU}$ algorithm with (Top) resilience index $s_1 = 0.201$ subject to attacks on $S_1 = 60$ agent sensors and (Bottom) resilience index $s_2 = 0.401$ subject to attacks on $S_2  = 120$ sensors.}\label{fig: MSIUResults}
\end{figure}

We examine the algorithm's performance with both fixed and time-varying attack sets, which are chosen uniformly at random. 
The adversary changes the measurement of each agent under attack ($n \in \mathcal{A}_t$) to $y_n(t) = -\theta^*$. The resilience of $\mathcal{SIU}$ does not depend on the adversary's strategy. 
Figure~\ref{fig: MSIUResults} shows that, for both fixed and time varying attack sets, and for both attacks on $S_1 = 60$ and $S_2 = 120$ agents, the $\mathcal{SIU}$ algorithm ensures that all of the agents' local estimates converge to $\theta^*$. 
The estimates converge more quickly when the resilience index is lower and there are fewer agents under attack (i.e., $S_1 = 60$). In general, there is a trade off between resilience and the rate at which the local estimates converge to the parameter of interest.

Recall, from the proof of Theorem~\ref{thm: mtsMain}, that the threshold $\gamma_t$ is an upper bound on the $\ell_2$ norm of the local estimation error, $\left\lVert x_n(t) - \theta^* \right\rVert_2$. We demonstrate the trade off between the algorithm's resilience index and the error bound $\gamma_t$. 
\begin{figure}[h!]
\centering
\includegraphics[width = \columnwidth]{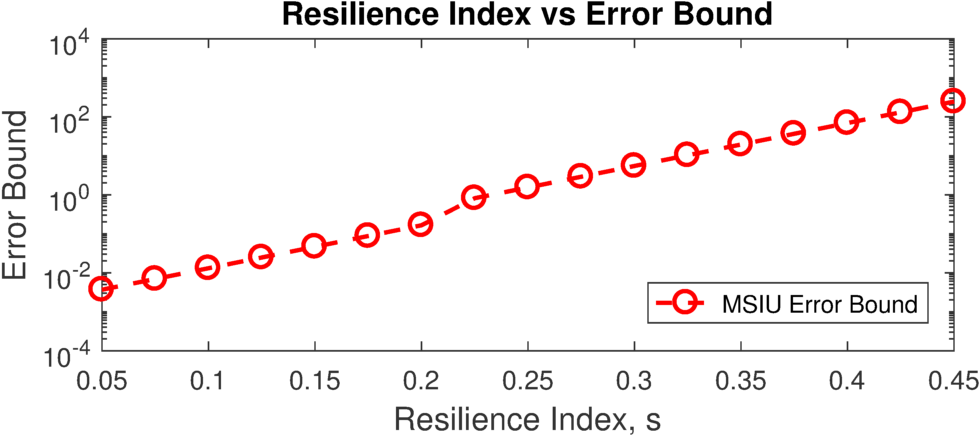}
\caption{Trade off between resilience index $s$ and the error bound $\gamma_t$ after $t = 500,000$ iterations of $\mathcal{SIU}$.}\label{fig: Resilience}
\end{figure}
Figure~\ref{fig: Resilience} shows that, as we increase the resilience index $s$, the error bound $\gamma_t$ (after $t = 500,000$ iterations) also increases. According to Theorem~\ref{thm: mtsMain}, the $\mathcal{SIU}$ algorithm ensures that the local estimates converge to $\theta^*$ eventually as long as less than $sN$ agents are under attack. In practice however, we may not be able to run an arbitrarily large number of iterations of $\mathcal{SIU}$, so we are interested in the error performance after a finite number of iterations.

\section{Conclusion}\label{sect: conclusion}
In this paper, we have presented the Saturated Innovation Update ($\mathcal{SIU}$) algorithm, a \textit{consensus + innovations} iterative distributed algorithm for resilient parameter estimation under sensor attacks. Under $\mathcal{SIU}$ for any (connected) network topology, local estimates of all agents converge polynomially to the parameter of interest if less than half of the agents' sensors are under attack. The number of tolerable attacks scales linearly with the number of agents, irrespective of the (connected) network topology. The $\mathcal{SIU}$ algorithm achieves the same level of resilience in a fully distributed setting as the most resilient centralized estimator. We demonstrated the performance of the $\mathcal{SIU}$ algorithms with numerical examples, and we showed that there exists a trade off between resilience and error performance. Future work will study resilient distributed estimation with more complicated sensing and communication models, e.g., communication networks with failing links, measurement noise, and time-varying dynamic parameters. 

\appendix
\subsection{Intermediate Results}
The proof of Lemma~\ref{lem: mi4} requires several intermediate results. The following results, from~\cite{Kar1} and \cite{Kar3}, study the convergence of scalar time-varying linear systems of the form
\begin{equation}\label{eqn: mtsSystem1}
	v_{t+1} = \left(1-r_2(t)\right) v_t + r_1(t),
\end{equation}
with $r_1(t)$, $r_2(t)$ of the form
\begin{equation}\label{eqn: mtsSystem1Rates}
	r_1(t) = \frac{c_1}{(t+1)^{\delta_1}},\: r_2(t) = \frac{c_2}{(t+1)^{\delta_2}},
\end{equation}
where $c_1, c_2 > 0$ and $0 < \delta_2 \leq  \delta_1 < 1$. 
\begin{lemma}[Lemma 25 in~\cite{Kar3}]\label{lem: mi0}
	Let $r_1(t)$ and $r_2(t)$ follow~\eqref{eqn: mtsSystem1Rates}. If $\delta_1 = \delta_2$, then there exists $B > 0$, such that, for sufficiently large non-negative integers, $j<t$, we have
	\begin{equation}\label{eqn: mi00}
		0 \leq \sum_{k = j}^{t-1} \left[ \left( \prod_{l = {k+1}}^{t-1}\left(1 - r_2(l) \right) \right) r_1(k) \right] \leq B,
	\end{equation}
	where the constant $B$ can be chosen independently of $t, j$. If $\delta_2 < \delta_1$, then, for arbitrary fixed $j$, we have
	\begin{equation}\label{eqn: mi01}
		\lim_{t \rightarrow \infty} \sum_{k = j}^{t-1} \left[ \left( \prod_{l = {k+1}}^{t-1}\left(1 - r_2(l) \right) \right) r_1(k) \right] = 0.
	\end{equation}
\end{lemma}
\noindent As a consequence of Lemma~\ref{lem: mi0}, for the system given by~\eqref{eqn: mtsSystem1}, if $\delta_1 = \delta_2$, then, $\left\lvert v_t \right \rvert$ remains bounded. If $\delta_1 > \delta_2$, then $v_t$ converges to $0$. 

The following result characterizes the rate of convergence of $v_t$ (given by~\eqref{eqn: mtsSystem1}) when $\delta_1 > \delta_2$. 
\begin{lemma}[Lemma 5 in~\cite{Kar1}]\label{lem: mtsIntermediate1}
	Consider the system in~\eqref{eqn: mtsSystem1} with $\delta_2 < \delta_1$. Then, we have
	\begin{equation}\label{eqn: mtsIntermediate1}
		\lim_{t\rightarrow\infty} (t+1)^{\delta_0} v_t = 0,
	\end{equation}
	for all $0 \leq \delta_0 < \delta_1 - \delta_2$ and for all initial conditions $v_0$.  
\end{lemma}

We present the following modification of Lemma~\ref{lem: mtsIntermediate1}.
\begin{lemma}\label{lem: mtsIntermediate2}
	Consider the scalar time-varying linear system:
	\begin{equation}\label{eqn: mtsSystem2}
		v_{t+1} = \left(1 - c_3 r_2(t) + c_4 r_1(t) \right) v_t + c_5 r_1(t),
	\end{equation}
	where $r_1(t), r_2(t)$ satisfy~\eqref{eqn: mtsSystem1Rates}, $\delta_1 > \delta_2$, and $c_3, c_4, c_5 >0$. The system in~\eqref{eqn: mtsSystem2} satisfies
	\begin{equation}\label{eqn: mtsIntermediate2}
		\lim_{t\rightarrow\infty} (t+1)^{\delta_0} v_t = 0,
	\end{equation}
	for all $0 \leq \delta_0 < \delta_1 -\delta_2$ and for all initial conditions $v_0$.
\end{lemma}

\begin{IEEEproof}
	Consider the expression $c_3r_2(t) - c_4 r_1(t)$. Since $\delta_1 > \delta_2$, for any $0 < \epsilon < \delta_1 - \delta_2$, there exists $T_0^{\epsilon}$, such that, for all $t \geq T_0^{\epsilon}$, we have 
	$c_3 r_2(t) - c_4 r_1(t) \geq (t+1)^{-\left( \delta_2 + \epsilon\right)}.$
Moreover, since $r_2(t)$ decreases in $t$ and $c_4 r_1(t) \geq 0$, there exists $T_1$ such that, for all $t \geq T_1$, we have
	$c_3 r_2(t) - c_4 r_1(t) \leq 1.$
Let $T = \max\left(T_0^{\epsilon}, T_1 \right)$. Thus, for all $t \geq T$, we have
\begin{equation}\label{eqn: mi2Proof3}
	(t+1)^{-\left( \delta_2 + \epsilon\right)} \leq c_3 r_2(t) - c_4 r_1(t) \leq 1,
\end{equation}
which means that, for all $t \geq T$, we have
\begin{equation}\label{eqn: mi2Proof4}
	\left\lvert v_{t+1} \right\rvert \leq  \left(1 - c_3 r_2(t) + c_4 r_1(t) \right) \left\lvert v_t \right\rvert + c_5 r_1(t).
\end{equation}
Let $w_T = \left\lvert v_T \right\rvert$, and, for $t \geq T$, define the time-varying linear system
\begin{equation}\label{eqn: mi2Proof5}
	w_{t+1} = \left( 1 - \left(t+1\right)^{-(\delta_2 + \epsilon)}\right)w_t + c_5 r_1(t). 
\end{equation}
As a consequence of~\eqref{eqn: mi2Proof3}, we have $\left\lvert v_t \right \rvert \leq w_t$ for all $t \geq T$. The system in~\eqref{eqn: mi2Proof5} falls under the purview of Lemma~\ref{lem: mtsIntermediate1}, which yields~\eqref{eqn: mtsIntermediate2}. 
\end{IEEEproof}

We now consider the system defined in~\eqref{eqn: miStackedSystem}:
\begin{equation*}
	\begin{split}
		v_{t+1} &= \left(1 - c_3 r_1(t) \right) v_t + c_4 r_1(t) w_t, \\
		w_{t+1} &= \left(1 - c_5 r_2(t) + c_6 r_1(t) \right) w_t + c_7 r_1(t) v_t,
	\end{split}
\end{equation*}
The following lemma shows that $v_t$ and $w_t$ remain bounded.

\begin{lemma}\label{lem: mi3}
	The system in~\eqref{eqn: miStackedSystem} satisfies
	\begin{align}
		\sup_{t \geq 0} \left\lvert v_t \right\rvert &< \infty, \label{eqn: mi3a} \\
		\sup_{t \geq 0} \left\lvert w_t \right\rvert &< \infty. \label{eqn: mi3}
	\end{align}
\end{lemma}

\begin{IEEEproof}
	\underline{Step 1}: We first show that $\sup_{t \geq 0} \left\lvert w_t \right\rvert < \infty$. Since $r_1(t)$ and $r_2(t)$ decrease in $t$, and, since $\delta_1 > \delta_2$, there exists (finite) $T_0 \geq 0$ such that, for all $t > T_0$, 
\begin{equation}\label{eqn: mi31}
	0 \leq 1 - c_3r_1(t) \leq 1,
\end{equation}
\begin{equation}\label{eqn: mi32}
	0 \leq 1 - c_5 r_2(t) + c_6 r_1(t) \leq 1.
\end{equation}
For all $t \geq T_0$, we can express $v_t$ as
\begin{equation}\label{eqn: mi33}
	\begin{split}
		v_t &= \prod_{j = T_0}^{t-1} \left(1 - c_3 r_1(j) \right) v_{T_0} + \\
	&\quad \sum_{j=T_0}^{t-1} \left[ \prod_{k = j+1}^{t-1}\left(1 - c_3r_1(k)\right) \right] c_4 r_1(j) w_j,
	\end{split}
\end{equation}
which means that, as a consequence of~\eqref{eqn: mi00},
\begin{equation}\label{eqn: mi34}
	\begin{split}
		\left\lvert v_t \right \rvert &\leq \left\lvert v_{T_0} \right \rvert + c_8 \sup_{l \in \left[T_0, t\right]} \left\lvert w_l \right \rvert,
	\end{split}
\end{equation}
for some constant $c_8 > 0$. 

\underline{Step 2}: From~\eqref{eqn: miStackedSystem}, we have
\begin{align}
\begin{split}\label{eqn: mi35}
	&\left\lvert w_{t+1} \right\rvert \leq \left(1 - c_5 r_2(t) + c_6 r_1(t) \right) \left\lvert w_t \right\rvert + c_7 r_1(t) \left\lvert v_t \right \rvert,
\end{split}\\
\begin{split}\label{eqn: mi37}
	&\quad \leq \left(1 - c_5 r_2(t) + c_9 r_1(t) \right) \sup_{l \in \left[T_0, t\right]} \left\lvert w_l \right\rvert + c_{10} r_1(t),\!\!
\end{split}
\end{align}
where $c_9 = c_6 + c_7$ and $c_{10} = c_7 \left\lvert v_{T_0} \right\rvert$. Define the system
\begin{equation}\label{eqn: mi38}
\begin{split}
	&m_{t+1}  = \\
	& \max \left(m_t, \left(1 - c_5 r_2(t) + c_9 r_1(t)\right) m_t + c_{10} r_1(t) \right),
\end{split}
\end{equation}
for $t \geq T_0$ with initial condition $m_{T_0} = \left\lvert w_{T_0} \right \rvert$. By definition, we have $m_t \geq \sup_{l \in \left[T_0, t \right]} \left\lvert w_l \right \rvert$. Define the system
\begin{equation}\label{eqn: mi39}
	\widetilde{m}_{t+1} = \left(1 - c_5 r_2(t) + c_9 r_1(t)\right) \widetilde{m}_t + c_{10} r_1(t),
\end{equation}
for $t \geq T_0$ with initial condition $\widetilde{m}_{T_0} = m_{T_0} = \left\lvert w_{T_0} \right\rvert$. Note that $\widetilde{m}_{T_0} \geq 0$ and that, since $t > T_0$, $1 - c_5 r_2(t) + c_9 r_1(t) \geq 0$, $\widetilde{m}_t$ nonnegative for all $t \geq T_0$. Further, note that system~\eqref{eqn: mi39} falls under the purview of Lemma~\ref{lem: mtsIntermediate2}, so we have
\begin{equation}\label{eqn: mi310}
	\lim_{t \rightarrow \infty} \widetilde{m}_t = 0.
\end{equation}

\underline{Step 3}: Since $\widetilde{m}_t$ is a nonnegative sequence that converges to $0$, there exists a time $T_1 \geq T_0$ such that $\widetilde{m}_{T_1 + 1} \leq \widetilde{m}_{T_1}$. Consider the smallest such choice of $T_1 \geq T_0$. 
By definition of $T_1$, we have
	$\widetilde{m}_{T_0} < \widetilde{m}_{T_0+1} < \dots < \widetilde{m}_{T_1}.$
Then, from~\eqref{eqn: mi38}, we have
	$m_t = \widetilde{m}_t$
for all $t \in \left[T_0, T_1\right]$. Thus, we have $m_t \leq m_{T_1}$ for all $t \in \left[T_0, T_1\right]$. Further, by definition of $T_1$ and $m_t$, we have
	$m_{T_1+1} = m_{T_1}.$
We now show that $m_{t} = m_{T_1}$ for all $t \geq {T_1}$. 

Define, for $t \geq T_1$,
\begin{equation}\label{eqn: mi314}
	J_1(t) = m_{T_1} - \left(1 - c_5 r_2(t) + c_9 r_1(t)\right) m_{T_1} - c_{10} r_1(t).
\end{equation}
We define $J_1(t)$ to be the difference between $m_{T_1}$ and $\left(1 - c_5 r_2(t) + c_9 r_1(t)\right) m_{T_1} + c_{10} r_1(t)$. If $J_1(t) \geq 0$ for all $t \geq T_1$, then, by~\eqref{eqn: mi38}, we have $m_t = m_{T_1}$ for all $t \geq T_1$.
By definition $m_{T_1} > 0$, so, after algebraic manipulation, we have
\begin{equation}\label{eqn: mi315}
	J_1(t) = \left( \frac{c_{11}}{(t+1)^{\delta_2}} - \frac{c_{12}}{(t+1)^{\delta_1}} \right) m_{T_1},
\end{equation}
where $c_{11} = c_5 c_2 > 0$ and $c_{12} = \left(c_6 + \frac{c_{10}}{m_{T_1}} \right) c_1 > 0$. Since $m_{T_1 + 1} = m_{T_1}$, $J(T_1)$ must be nonnegative, which is true if and only if
\begin{equation}\label{eqn: mi316}
	T_1 \geq \left(\frac{c_{12}}{c_{11}}\right)^{\frac{1}{\delta_1 - \delta_2}} - 1,
\end{equation}
so all $t \geq T_1$ are also at least the right hand side of~\eqref{eqn: mi316}, which means that $J_1(t) \geq 0$ for all $t \geq T_1$. 

\underline{Step 4}: Thus we have $m_{t} = m_{T_1}$ for all $t \geq T_1$. Using the fact that $m_t \leq m_{T_1}$ for all $t \in \left[T_0, T_1 \right]$, and $m_t \geq \sup_{l \in \left[T_0, t\right]} \left \lvert w_l \right \rvert$, we have
\begin{equation}\label{eqn: mi317}
	\sup_{t \geq T_0} \left\lvert w_t \right\rvert \leq m_{T_1} < \infty.
\end{equation}
Since $T_0 < \infty$, we also have
\begin{equation}\label{eqn: mi318}
	\sup_{t \in \left[0, T_0\right]}\left\lvert w_t \right\rvert < \infty.
\end{equation}
Combining~\eqref{eqn: mi317} and~\eqref{eqn: mi318} yields~\eqref{eqn: mi3}. Now, we show that $\sup_{t \geq 0} \left\lvert v_t \right\rvert < \infty$. Let $\sup_{t \geq 0} \left\lvert w_t \right \rvert = B_w < \infty$. By~\eqref{eqn: mi34}, we have
\begin{equation}\label{eqn: mi319}
	\sup_{t \geq T_0} \left\lvert v_t \right \rvert \leq \left\lvert v_{T_0} \right \rvert + c_8 B_w < \infty.
\end{equation}
Since $T_0 < \infty$, we also have
\begin{equation}\label{eqn: mi320}
	\sup_{t \in \left[0, T_0\right]}\left\lvert v_t \right\rvert < \infty.
\end{equation}
Combining~\eqref{eqn: mi319} and~\eqref{eqn: mi320} yields~\eqref{eqn: mi3a}.
\end{IEEEproof}

\subsection{Proof of Lemma~\ref{lem: mi4}}

\begin{IEEEproof}
	We first prove~\eqref{eqn: mi4a}, that $\lim_{t \rightarrow \infty} (t+1)^{\delta_0} w_t = 0$ for all $0 \leq \delta_0 < \delta_1 - \delta_2$. From Lemma~\ref{lem: mi3}, we have $\left\lvert v_t \right \rvert \leq B_v < \infty$. Then, for sufficiently large $t$, we have
\begin{equation}\label{eqn: mi41}
	\left\lvert w_{t+1} \right \rvert \leq \left(1 - c_5 r_2(t) + c_6 r_1(t) \right) \left\lvert w_t \right \rvert + c_7 B_v r_1(t).
\end{equation}
The recursion in~\eqref{eqn: mi41} falls under the purview of Lemma~\ref{lem: mtsIntermediate2}, and~\eqref{eqn: mi4a} immediately follows. Further, as a consequence of Lemma~\ref{lem: mtsIntermediate2}, there exists $R_w > 0$ such that
	$\left\lvert w_t \right \rvert < R_w (t+1)^{\delta_0},$
for all $0 \leq \delta_0 < \delta_1-\delta_2$. We now prove~\eqref{eqn: mi4}. Since $\left\lvert w_t \right \rvert < R_w (t+1)^{\delta_0}$, we have, for sufficiently large $t$,
\begin{equation}\label{eqn: mi43}
	\left \lvert v_{t+1} \right \rvert \leq \left(1 - c_3 r_1(t) \right) \left\lvert v_t \right \rvert + \frac{c_1 c_4 R_w}{(t+1)^{\delta_1 + \delta_0}}.
\end{equation}
The recursion in~\eqref{eqn: mi43} falls under the purview of Lemma~\ref{lem: mtsIntermediate1}, and we have
\begin{equation}\label{eqn: mi44}
	\lim_{t\rightarrow\infty} (t+1)^{\delta_0'} v_t = 0,
\end{equation}
for all $0 \leq \delta_0' < \delta_0$. Taking $\delta_0$ arbitrarily close to $\delta_1 - \delta_2$ yields~\eqref{eqn: mi4}. 
\end{IEEEproof}

\bibliography{IEEEabrv,References}

\end{document}